\input amstex

\input amsppt.sty
\input epsf

\hsize=5.00in
\vsize=8.10in
\parindent 20 pt
\NoBlackBoxes

\define \be{\beta}
\define \Dl{\Delta}
\define \dl{\delta}
\define \g{\gamma}
\define \G{\Gamma}

\define \ve{\varepsilon}
\define \vp{\varphi}

\define \si{\sigma}

\define \df{\dsize\frac}

\define \fa{\forall}
\define \fc{\frac}

\define \la{\langle}
\define \ra{\rangle}
\define \p{\partial}

\define \ri{\rightarrow}

\define \ov{\overline}

\define \edm{\enddemo}
\define \ep{\endproclaim}

\define \bk{\bigskip}

\define \1{^{-1}}
\define \2{^{-2}}

\define \CP{\Bbb C\Bbb P}

\define \CPt{\Bbb C\Bbb P^2}
\define \BC{\Bbb C}

\define \BP{\Bbb P}

\define \BZ{\Bbb Z}

\define \tX{\tilde{X}}
\define \tY{\tilde{Y}}

\define \Cen{\operatorname{Center}}

\baselineskip 20pt

\topmatter
 \rightheadtext{Braid groups, algebraic surfaces and fundamental groups}
\title Braid groups, algebraic surfaces and fundamental groups of complements of branch
 curves\endtitle
\author M. Teicher\endauthor
\abstract\nofrills{\bf Abstract.} An overview of the braid group techniques in the
theory of algebraic surfaces from Zariski to the latest results is presented.
An outline of the Van Kampen algorithm for computing fundamental groups of
complements of curves and  the modification of Moishezon-Teicher regarding branch curves
of generic projections  are given.
The paper also contains a description of a quotient of the braid group, namely $\tilde
B_n$ which
  plays an important role in the description of fundamental groups of
complements of branch curves.
It turns out that all such groups are ``almost solvable'' $\tilde B_n$-groups.
Finally, the possible applications to study moduli spaces of surfaces of general type
are described and new examples of positive signature spin surfaces whose fundamental
groups can be computed using the above algorithm (Galois cover of Hirzebruch surfaces)
are presented.\endabstract 
 \thanks This research was partially supported by the Emmy
Noether Mathematics Research Institute, Bar-Ilan University, Israel.\endthanks
\endtopmatter \baselineskip 20pt \subheading{0. Introduction}

This manuscript is based on our talk in Santa Cruz, July 1995.
It presents the applications of the braid group technique to the study
of algebraic surfaces  and curves in general and to the moduli space of
surfaces and the topology of complements of curves in
particular. 
These techniques started with Enriques, Zariski and Van Kampen in the 30's
(see \cite{VK}, \cite{Z}) and were revived by Moishezon in the late 70's (see,
e.g., \cite{Mo1}).
The manuscript  includes a survey on the topology of complements of branch
curve starting with Zariski's results, as well as new results (related to a quotient $\tilde B_n$ of 
the braid group)  and   an open
question on the topic.   The manuscript is divided as follows:

\roster \item"0."  Introduction
\item"I."  The connections between  classification of
algebraic surfaces and related fundamental groups
\item"II." Known   results of fundamental
groups of complements of branch curves; an open question
 \item"III." Presentation of $\tilde B_n$, a quotient
 of the braid group 
 \item"IV." Two new theorems on the fundamental groups of  complements of branch curves of
$V_3$ (Veronese of order 3)
 \item"V." An algorithm to compute
fundamental groups of complements of   branch curves
 \item"VI." The Braid monodromy (Step (b) of the algorithm)
\item"VII." The Enriques-Van Kampen method  (Step (c) of the algorithm)
\item"VIII." Some facts on the structure of $\tilde B_n$ and $\tilde B_n$-groups
(steps (e) and (f) of the algorithm)
 \item"IX." The connection between fundamental groups
of complements of branch curves and Galois covers
 \item"X."  Galois covers of Hirzebruch surfaces: new examples\endroster 

\bigskip

\subheading{I.  The connections between classification of
algebraic surfaces and related fundamental groups}

In 1977 Gieseker prove that the  moduli space of surfacaes of general type is a quasi-projective variety (see
\cite{G}).
Unlike the case for curves it is not irreducible.
 Catanese and Manetti proved results about the structure
and the number of components of moduli spaces (see, e.g., 
\cite{C1}, \cite{C2}, \cite{C3}, \cite{C4}, \cite{C5}, \cite{C6}, \cite{CCiLo},
\cite{CW}, \cite{Ma}). 
Not much is known about these moduli spaces.
Nevertheless, unlike previous expectations, simply connected (and spin) surfaces exist
also in the
 $\tau>0$ area, $\tau=\frac{1}{3}(C_1^2-2C_2)$ (see
\cite{MoTe1}, \cite{MoTe2}, \cite{MoTe3}, \cite{Ch}, \cite{MoRoTe},
\cite{PPX}).

The fact that algebraic surfaces are nontrivial geometric objects was
remarkably confirmed  by S. Donaldson who showed that among algebraic surfaces
one can find homeomorphic non-diffeomorphic (simply-connected)  4-manifolds.  
In particular, he produced the first counterexamples to the h-cobordism
conjecture in dimension four.  
Donaldson's theory was also used  to construct the first examples of homeomorphic
non-diffeomorphic (simply-connected) algebraic surfaces of general type 
(\cite{FMoM}, \cite{Mo2}).
In 1994, Witten \cite{W} and later Witten and Sieberg \cite{SW} defined a new set of
invariants for $4$-manifolds (monopole invariants), and have shown the equivalence of
this invariant with Donaldson's polynomial.
These invariants take a simple form for K\"ahler surfaces.

	We expect that the connected components of moduli spaces of algebraic
surfaces  (of general type) correspond to the principal diffeomorphism classes
of corresponding topological 4-manifolds.  
Thus, it is possible that Donaldson's polynomial invariants will distinguish
these connected components.  
However, we present here  a more direct geometrical approach.

 \bk The ultimate goals of  the braid group
techniques are 
finding  new invariants distinguishing connected components of the moduli
space of surfaces of general type.  For that we try to
  compute different fundamental groups related to the surface, groups which do
not change when one moves in a connected component of the moduli space.
The first groups we
compute are $\pi_1(\BC^2-S)$ and $\pi_1(\CPt-S)$ where $S$ is the branch
curve of a generic projection $X\to\CPt.$
If $\pi_1$ is ``big''  then it can distinguish between connected
components. 
If they are ``small'' there is hope to compute $\pi_2$ as a module over
$\pi_1.$ 
We can also compute fundamental groups of surfaces of general type.
This is especially interesting in the  positive signature   area which is still rather
wild.

For  minimal surfaces of general type it
turns out that all the information is contained in the canonical class: i.e. it
is a diffeomorphism invariant and all other information about Donaldson's
polynomials must follow from it. Thus,
for the problem of finding invariants of
deformation types of surfaces of general type   we are almost 
where we were 15 years ago (the only new invariant is divisibility of the canonical 
class). So fundamental groups of the complements to branching curves of
generic projections might still be the best bet for this subject.

We want to recall here that computing fundamental groups of complements of a
plane curve is enough in order to  understand the topology of a complement
in $\Bbb P^N$ of any algebraic subset (as proven by Zariski).
In fact, for a generic $\Bbb P^2 $ in $\Bbb P^N:$
$$ \pi_1(\Bbb P^N-V)\simeq\pi_1(\Bbb P^2-\Bbb P^2\cap V).$$
Furthermore, we recall that lately there is also a growing interest in fundamental
groups of algebraic varieties in general.
A very partial list includes \cite{BoKa}, \cite{CMan}, \cite{DOZa} \cite{L1}, \cite{L2},
\cite{Si}, \cite{To}.

The braid group appears in the formulation of the results and as an essential step of
the algorithm for computing   fundamental groups   of complements of curves (see Section
V). 

\bk

\subheading{II. Known results on fundamental groups of complements of branch
curves; an open question}

Consider the following situation:
$$\align\text{Surface}\ &X\hookrightarrow \BC\BP^N\\
&\downarrow\text{generic projection}\\
\quad S\subseteq\      &\CPt\qquad S= \text{branch curve}\endalign$$ We
denote: $G=\pi_1(\BC^2-S,*),$\quad $\ov G=\pi_1(\CPt-S,*).$

We want to  to find a general formula for $G$ and $\ov G$
which depends on known invariants of $X.$ 
As we said in the our introduction, the topic started with Zariski who proved
in the 30's that if  $X$ is a cubic surface in $\BC\BP^3$ then
$\ov G\simeq Z_2\star Z_3$ (see \cite{Z}).
In the late 70's Moishezon proved that if $X$ is a $\deg n$ surface in
$\BC\BP^3$ then $G\simeq B_n,$\quad $\ov G\simeq B_n/\Cen$ (see \cite{Mo1}).
In fact, Moishezon's result for  $n=3$ is the same as
Zariski's result since  $B_3/\Cen\simeq Z_2\star Z_3.$

The next example was $V_2$ (Veronese of order 2) (see \cite{MoTe3}).
In all the above examples we have $G\supset F_2$ where $F_2$ is a free noncommutative subgroup
with 2 elements.
We call a group $G$ ``big'' if $G\supset F_2.$

Since 1991  the following examples have been discovered:
$V_3$,   the Veronese of order 3
which was done by Moishezon and Teicher in
\cite{MoTe7}, \cite{MoTe8}, \cite{MoTe9}, \cite{MoTe10}, \cite{MoTe11}, \cite{Te2}, and
generalized later to general $V_n$ (preprint);
$X_{ab},$  the embedding of $\BC\BP^1\times \BC\BP^1$ into $\BC\BP^N$ w.r.t.
a linear system $|a\ell_1+b\ell_2|$; \ $CI$, the complete intersection which
was done by A. Robb in his Ph.D. Thesis in 1994, (see \cite{Ro}).

\bigskip
Unlike previous expectations, in all the new example    $G$ is not ``big''.
Moreover, $G$ is ``small'', i.e.,
$G$ is ``almost solvable'', i.e., it contains a subgroup of finite index
which is solvable. 
It turned out that there exists a quotient of the braid group (by a subgroup ofthe
commutant), namely $\tilde B_n$ s.t.
 all new results give $ G=\tilde
B_n$-group and $\ov G=G/$central element $(\tilde B_n$-group is a group on which $\tilde B_n$
act).  
For $CI,$ \ $G$ is $\tilde B_n$ itself.
So the old examples were exceptions ($V_2$ often turns out  to be an exception) and 
fundamental groups of complements
of branch curves are not ``big''.  They are surprisingly ``small''. 
Moreover, in all the new examples $G,\ov G$    are an extension of a solvable group by a symmetric one.
Based on that fact we ask the following
\proclaim{Question}
For which familes of simply connected  algebraic surfaces of
general type is the fundamental group of the complement of the branch curve of
a generic projection to  $\CPt$  an extension of a solvable group by a
symmetric group?\endproclaim

We believe that the answer to this question lies in the decomposable structure of the corresponding
4-manifold.
One should also notice that if a group $G$ is ``big'' then it is not ``small''
 and if it is ``small''
then it is not ``big''.

\bk

\subheading{ III. Presentation of $\tilde B_n$, a quotient
 of the braid group }

The braid group is connected to fundamental groups of complements of branch
curves in two ways.  The first way is through the appearance of its quotient $\tilde
B_n$ in the description of such groups (see  Section II), and the second way is
through the use of the braid group as a major tool in the algorithm for computing such
groups (see Section IV).

We first review the definition of  braid group (see also \cite{A} and \cite{B}), and then
we shall define its quotient $\tilde B_n.$
We will work with a geometric model of the braid group.

\newpage

\definition{Definition:\ The braid group
 $B_n$}

Let $D$ be a topological disc, $K\subset D$ finite.
Consider:
\ $\{\be|\be: D\to D$ diffeomorphism, \ $\be(K)=K,$\ $\be|_{\p D}=Id\}.$
Clearly, $\{\be\} $ is a group which acts naturally on $\pi_1(D-K).$
We define an equivalence relation on $\{\be\}$ as follows:
$\be_1\sim \be_2\Leftrightarrow$ the action of $ \be_1,\be_2$   on $\pi_1(D-K)$ coincide.
 $B_n=\{\be\}\bigm/_\sim$

We have to distinguish   certain elements in $B_n.$\enddefinition

\definition{Definition:\ Half-twist w.r.t.
$\left[-\df{1}{2},\frac{1}{2}\right]$}

Consider $D_1,$ the  unit disc, $\pm\df{1}{2}\in D_1.$
Take $\rho: [0,1]\to[0,1]$ continuous s.t. $\rho(r)=\pi$\quad
$r\le\df{1}{2}$\quad $\rho(1)=0.$
Define $\dl: D_1\to D_1: \dl(re^{i\theta})=re^i(\theta+\rho(r)).$
Clearly, $\dl\left(\df{1}{2}\right)=-\df{1}{2},$\quad
$\dl\left(-\df{1}{2}\right)=\df{1}{2},$ and $\dl|_{\p D_{1}}=Id.$
The disc of radius $\df{1}{2}$ rotates $180^\circ$ counterclockwise.
Outside of this disc it rotates in smaller and smaller angles till it rests
on the unit circle.
Thus we get a braid $[\dl]\in
B_2\left[D_1,\left\{\pm\df{1}{2}\right\}\right].$ $[\dl]$ is called the
half-twist w.r.t. the segment
$\left[-\df{1}{2},\df{1}{2}\right].$\enddefinition \bigskip
Using the above definition we define a generalized
half-twist.
\definition{Definition: $H(\si),$  half-twist w.r.t. a path $\sigma$}

\flushpar Let $D,\ K$ be as above, $a,b\in K.$
Let $\sigma$ be a path from $a$ to $b$ which does not meet any other
point of $K.$   We take $D_2$ a small topological disc in $D$ s.t.
$\sigma\subset D_2\subset D,\quad D_2\cap K=\{a,b\}.$
We take $\psi: D_2\to D_1$ (unit disc) s.t. 
$\psi(\sigma)=\left[-\fc{1}{2},\fc{1}{2}\right].$
\ $\psi(a)=-\fc{1}{2}$\quad $\psi(b)=\fc{1}{2}.$
 We consider a ``rotation'' $\psi\dl\psi\1: D_2\to D_2;$
 $\psi\dl\psi\1$ is identity on the boundary of $D_2.$
We extend it to $D$ by identity.
 $H(\sigma)=$ [extension of $\psi\dl\psi\1].$
\enddefinition

\centerline{
\epsfysize=1.1in
\epsfbox{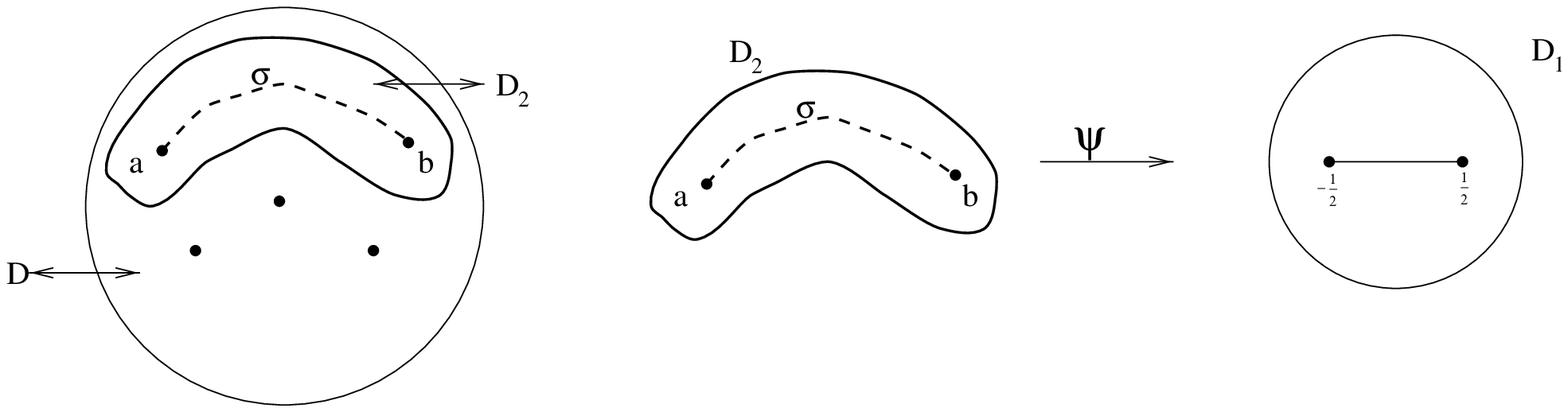}}  

\midspace{.05in}

We shall present now  $\tilde B_n,$ the quotient of the braid group by
commutators of the transversal half-twists.
We define: 
\definition{Definition: Transversal half-twists}
\flushpar  $H(\si_1)$ and $H(\si_2)$ are transversal if
$\si_1\cap\si_2=\{$one point which is not an end point$\}$

\medskip

\centerline{
\epsfysize=.65in
\epsfbox{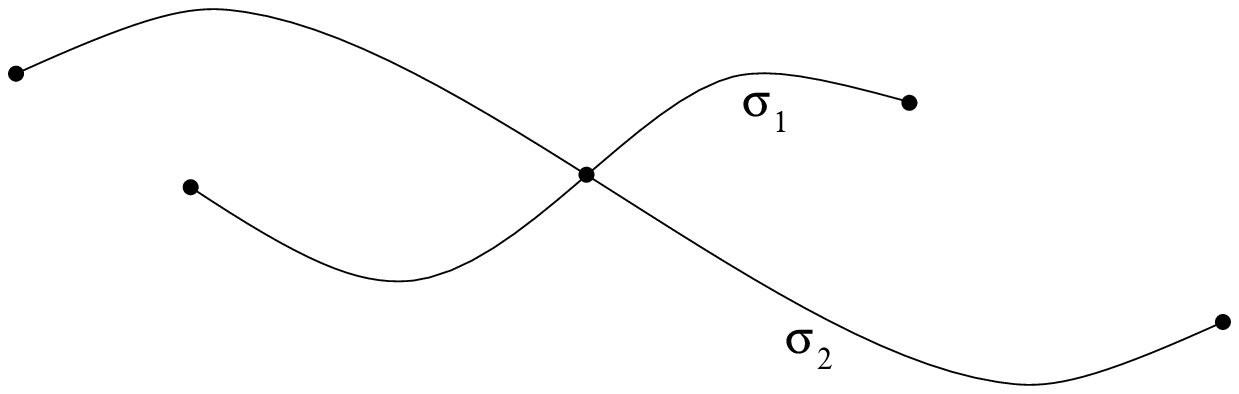}}

\enddefinition

\definition{Definition: $\tilde B_n$}
\flushpar  Let $X,Y$ be a pair of transversal half-twists.
   Let $[X,Y]=XYX\1Y\1.$
Let $\la [X,Y]\ra$ be the subgroup normally generated by $[X,Y].$
   $\tilde B_n=B_n\bigm/_{\la[X,Y]\ra}.$\enddefinition

\remark{Remark}
Since all transversal half-twists are conjugated,
  $\la[X,Y]\ra$  contains every commutator of transversal half-twists and thus
 $\tilde B_n$ is independent of the choice of
$X,Y.$\endremark
 One can find a description of $\tilde B_n$ in \cite{Te1}.
\bk

 \subheading{IV. Two new theorems on the fundamental groups of  complements of branch curves of
$V_3$ (Veronese of order 3)}

For example, we shall formulate exactly the structure theorem concerning
$V_3.$

\proclaim{Theorem 1} \rom{\cite{MoTe9}}

\flushpar Let $X$ be $V_3$ (the Veronese of order $3).$

\flushpar Let $S$ be the branch curve of a generic projection to $\CPt.$
Then:

$$\pi_1(\BC^2-S,*)\simeq \tilde B_9\ltimes G_0(9)\bigm/_{N_9}$$
where
$$ G_0(9) = \text{central extension of a free group with $8$ elements}
=\la u_1,\dots,u_8,\tau\ra$$
s.t.
$$ [u_i,u_j]=\cases \tau\quad & |i-j|=1\\
1\quad & |i-j|\ne 1\endcases$$
$$\tau\in \Cen G_0(9),\quad \tau^2=1.$$
\flushpar There exists a standard base of $\tilde B_9:\tX_1,\dots,\tX_8$ s.t.
 the action of $\tilde B_9$ on $G_0(9)$ is as follows:

\flushpar$(u_i)_{x_{k}}=\cases u_i\tau\quad & k=i\\
u_i\quad & |i-k|\ge 2\\
u_ku_i\quad & |i-k|=1.\endcases$

\flushpar$ N_9=\la u_i^3=X_i^3,\quad \tau=c\ra$ where  $c\in \Cen\tilde
B_9,$\quad $c^2=1.$\ep \bigskip
The ``almost solvable'' theorem concerning $V_3$ is as follows:
\proclaim{Theorem 2} \rom{\cite{MoTe10}}
\flushpar Let $X,$\ $S$ be as in the previous theorem.
\flushpar Let $G=\pi_1(\BC^2-S).$ Then there exists a series $1<
H_{9,0}'<H_{9,0}<H_9<G$
s.t.
$$\align
& G/H_9\simeq S_9\\
&H_9/H_{9,0}\simeq \BZ\\
&H_{9,0}/H_{9,0}'\simeq (\BZ\oplus\BZ/3)^8
\quad\\
&H_{9,0}'(=H_9')=\{1,c\}\simeq \BZ/2\quad (c\in \Cen G).
\endalign$$
\ep
We did not discuss yet where does the braid group enter into the
calculation of fundamental groups of complements of branch curves; we do this in
the next section.

\bigskip

\subheading{V. An algorithm to compute
fundamental groups of complements of   branch curves}

In this section we state the main steps used so far for computing such
groups: \roster\item"(a)" Degeneration of the surface to a union of planes
where no 3 planes meet in in a line. \item"(b)" Computing the braid monodromy
of the branch curve (using the above degeneration).
\item"(c)" Enriques-Van Kampen method for getting a finite presentation of   $\pi_1(\CPt-S)$ (using the braid
monodromy).
\item"(d)" Invariance properties of the braid monodromy (to produce more
relations in\linebreak  $\pi_1(\CPt-S)$ than those induced from the Van Kampen method).
\item"(e)" Studying $G$ as a    $\tilde B_n$-groups and looking for prime elements.
\item"(f)" Proving ``almost solvability'' when available.
\endroster

At the moment we work on eliminating the condition that no 3 planes meet in a
line in order to enlarge the variety of surfaces to which we can apply our
methods. 
The reason that we need the degeneration at all is to simplify the
computations of the braid monodromy of the branch curve.
If the surface is degenerated to a union  of planes where no 3 planes meet in
a line, then the degenerated object has a branch curve which is partial to an
arrangement of lines known as ``dual to a generic''.

An arrangement of lines ``dual to generic'' is an arrangement in which there are exactly
2 multiple points (where $m$ lines meet, $m\ge 3)$ on every line.
In \cite{MoTe4} we presented an algorithm for computing the braid monodromy of
arrangement ``dual to a generic''
In \cite{MoTe6} we presented an algorithm how to get from the braid monodromy
of the degenerated braid curve, the braid monodromy of the original curve.
To eliminate the condition in (a) means to produce an algorithm for computing
braid monodromies of arrangements of lines which are not ``dual to generic''.
This as explained  earlier will enlarge the variety of surfaces for which we
can compute  $\pi_1(\CPt-S)$.

\bk

\subheading{VI. The braid monodromy (Step (b) of the algorithm)}

Computing the braid monodromy is the main tool to compute fundamental groups
of complements of curves (Step (b)).
In this section we define the braid monodromy and compute some examples.
\definition{Definition:\ The braid monodromy w.r.t. $S,\pi,u$}

\flushpar Let $S$ be a curve, $S\subseteq \Bbb C^2$

\flushpar Let $\pi: S\to\BC^1$ be defined by
 $\pi(x,y)=x.$
We denote $\deg\pi$ by $m.$
\flushpar Let $N=\{x\in\Bbb C^1\bigm| \#\pi\1(x)\lneqq m\}.$
   Take $u\notin N,$ s.t.  $x\ll u$ \ $\forall x\in N.$
Let  $ \BC^1_u=\{(u,y)\}.$
\flushpar  There is a  natural defined homomorphism $\pi_1(\BC^1-N,u)\overset
\vp\to\ri B_m[\BC_u^1,\BC_u^1\cap S]$ which is called {\it the braid
monodromy w.r.t.} $S,\pi,u.$\enddefinition

\medskip

\centerline{
\epsfysize=3in
\epsfbox{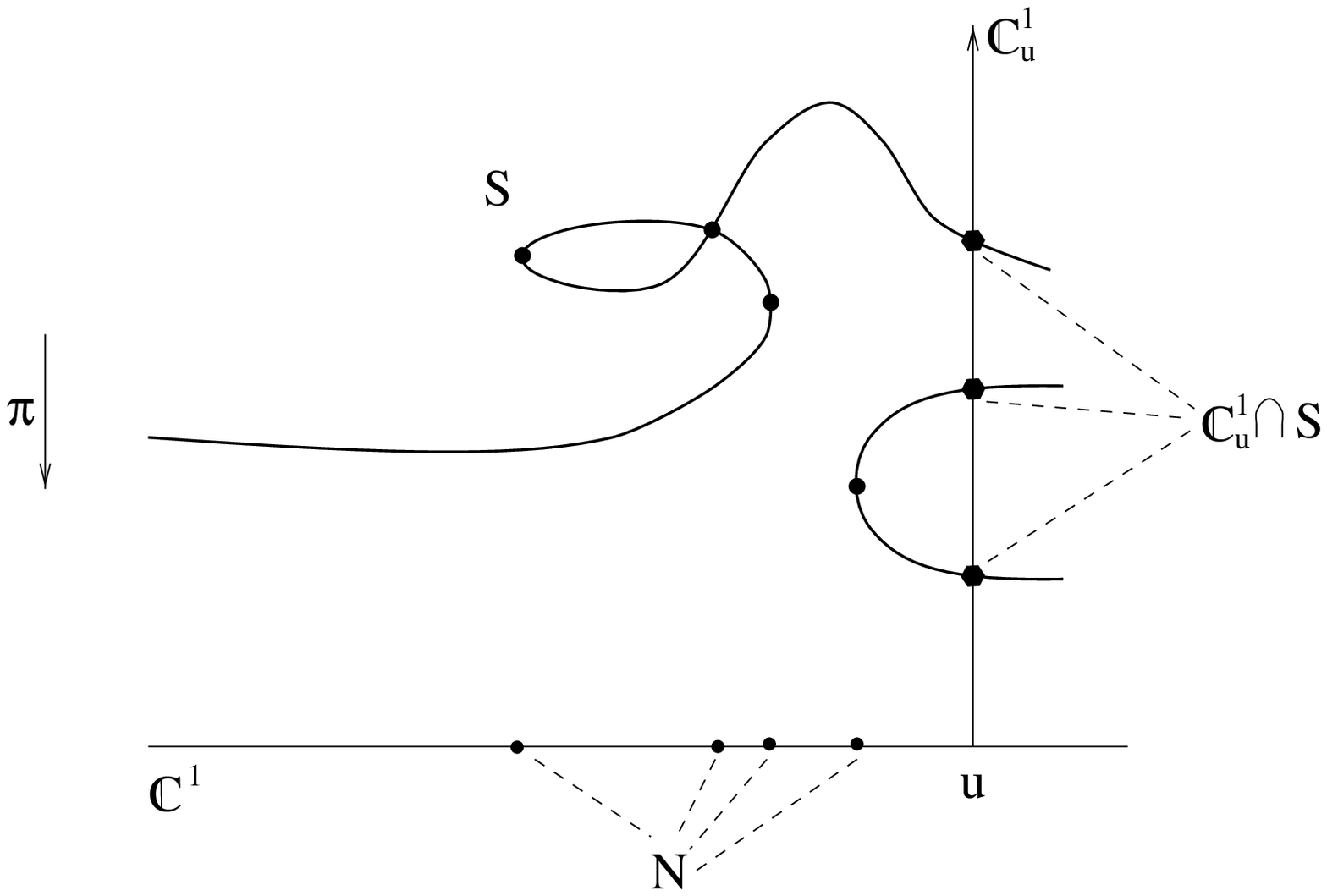}}  

\medskip

\remark{Remark}
The classical monodromy factors through the braid monodromy
$$\alignat2 \pi_1(\BC^1-&N,u)\to&&B_m[\BC_u^1,\BC_u^1\cap S]\\
&\searrow&&\downarrow\\
& && S_m\endalignat$$\endremark
 
\bk
\subheading{Example of computing braid monodromy of a curve with only one
singular point}

\flushpar Let  $S$ be defined by $ y^2=x^m.$
\flushpar For  $\pi: S\to \BC^1$ defined by $\pi(x,y)=x$ we have $\deg\pi=2. $
$S$ has only one singular point $(0,0)$ and thus
 $N=\{0\}.$ We take $u=1$.  
Clearly,  $\BC_u^1\cap\Bbb S=\{-1,1\}. $
\flushpar
Let   $\dl(t)=e^{2\pi it}$  
($\delta(t)$ is a closed loop that starts in $u).$
$\dl$ is a generator of \linebreak $\pi_1(\BC^1-N,u).$

\flushpar We lift $\delta(t)$ to $S.$
There are 2 liftings:

 $\dl_1(t)=\left(e^{2\pi it},e^{\fc{2\pi
itm}{2}}\right)$

  $\dl_2(t)=\left(e^{2\pi it},-e^{\fc{2\pi itm}{2}}\right)$
\flushpar The projections of $\delta_1(t)$ and $\delta_2(t)$ to $\BC_u^1$ are:

$a_1(t)=e^{\fc{2\pi itm}{2}}=(e^{\pi it})^m$

  $a_2(t)=-e^{\fc{2\pi itm}{2}}=- (e^{\pi it})^m$

\flushpar By definition of the braid monodromy, $\vp(\delta)$ is induced from
the motion \linebreak $\{(e^{\pi i t})^m,-(e^{\pi it})^m\}.$

 \flushpar
The braid  induced from the motion $\{e^{\pi it},-e^{\pi it}\}$ is $H=H([-1,1])$  which is
the half-twist in $\BC_u^1$ w.r.t. $[-1,1].$
Clearly, $\vp(\dl)=H^m.$

The above example is almost a proof for the following theorem of Zariski.

\proclaim{Theorem} \rom{(Zariski)}
Let $S$ be a cuspidal curve.
Assume that above each point of $N$ there is only one singular point of $\pi.$
Let $x_0\in N.$
Let $\dl$ be a loop in $\pi_1(\BC-N,u)$ around $x_0$ \ $(\dl$ is simple and
no other point of $N$ is inside $\dl).$
Let \newline $\vp:\pi_1(\BC-N,u)\to B_m$ be the braid monodromy.
Then $\vp(\dl)=H^\ve$ where $H$ is a half-twist and
$$\ve=\cases1\quad & (x_0,y_0)\ \text{is a branch point of}\ \pi\\
2\quad & (x_0,y_0)\ \text{is a node of}\ S\\
3\quad & (x_0,y_0)\ \text{is a cusp of}\ S\endcases$$\ep
\remark{Remark}
Clearly, the complexity in finding  $\vp(\dl)$ lies in finding
$H.$\endremark 

\bk

\subheading{VII. The Enriques-Van Kampen method (Step (c) of the algorithm)}

The Van Kampen method gives us a finite presentation in terms of generators and relation of 
plane complements of curves.

The categorical version of the Van-Kampen Theorem is as follows:

\proclaim{Van Kampen Theorem} \rom{\cite{VK}}
Let $\ov S\subseteq \CPt$ be a projective curve, which is transversal to the
line at infinity.
Let $S=\ov S\cap\BC^2.$
Let  $\vp_u:\pi_1(\BC-N,u)\to B_m[\BC_u^1,\BC_u^1\cap S]$ be the  braid monodromy w.r.t. 
$S,\pi,u.$
Then:
\roster\item"(a)" $\pi_1(\BC^2-S,*)=\pi_1(\BC_u^1-S,*)/\{\be(V)=V|\be\in\operatorname{Im} 
\vp_u,\ V\in
\pi_1(\BC_u^1-S)\}.$
\item"(b)" $\pi_1(\CPt-\ov S)\simeq
\pi_1(\BC^2-S)/\la\G\ra$ where $ \G $ is a simple loop in $\Bbb C_u^1-S$ around $S\cap 
\Bbb C_u^1=\{q_1,\dots,q_m\}.$\endroster\ep

We want to rephrase the Van Kampen theorem in a way that it can be used with greater facility.
To that end we need the notion of a good geometric
base for the fundamental group of a punctured disc.
We recall that for $D-K,$ a   punctured disk,  $\pi_1(D-K)$ is a free group and $B_n[D,K]$ acts naturally on
$\pi_1(D-K).$
Before defining a good geometric base we need 2 additional definitions.

 \definition{Definition:\ $\ell(q)$} 

\flushpar Let $D$ be a topological disc, $K\subset D,$\ $K$ finite,
$u\in\p D.$
\flushpar Let $a\in K,$ \ $q$ a simple path from $u$ to $a$ such that $q\cap
K=a.$ \flushpar  Let $c$ be a simple loop equal to the (oriented) boundary of a
small neighborhood $V$ of $a$ chosen such that $q'=q-V\cap q$ is a simple path.
Then $\ell(q)=q'\cup c\cup q'{}\1$ (see figure).

\medskip

\centerline{
\epsfysize=1.1in
\epsfbox{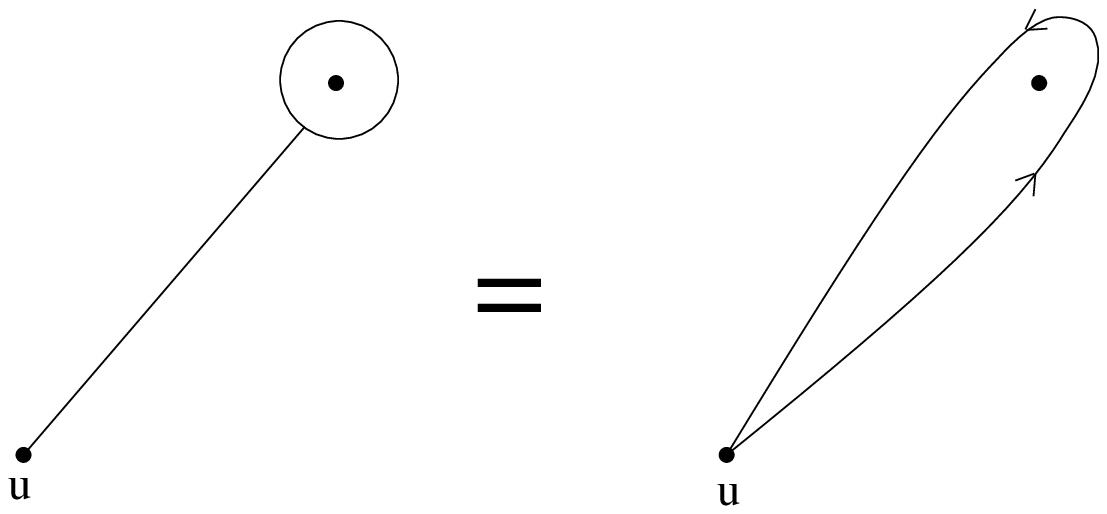}}  

\medskip

\flushpar  We use the same notation $\ell(q)$ also for the element of
$\pi_1(D-K,u)$ corresponding to $\ell(q).$\enddefinition

\definition{Definition:\ A bush}

\flushpar Let $D,K,u$ be as above.
Let $K=\{a_1,\dots, a_n\}.$
\flushpar  Consider in $D$ ordered sets of simple paths $(T_1,\dots, T_n)$
connecting $a_i$'s with $u$ such that
\roster\item $\fa i=1,\dots,n\ t_i\cap w_j=\emptyset$ if $i\ne j;$
\item $\operatornamewithlimits\cap\limits_{i=1}^n T_i=u;$
\item for a small circle $c(u) $ around $u$ each $u_i'=T_i\cap c(u)$ is a
single point and the order in $(u_1',\dots,u_n')$ is consistent with the
positive (``counterclockwise'') orientation of $c(u).$\endroster
We say that two such sets $(T_1,\dots,T_n)$ and $(T_1',\dots,T_n')$ are
equivalent if $\fa i=1,\dots,n$:
$\ell(T_i)=\ell(T_i')\quad \text{(in}\ \pi_1(D-K,u)).$
\flushpar An equivalence class of such sets is called a bush in $(D-K,u).$
The bush represented by $(T_1,\dots,T_n)$ is denoted by $\la
T_1,\dots,T_n\ra.$
\enddefinition

\definition{Definition:\ A good geometric base  ($g$-base)}

Let $D,K$ be as above. A good geometric base of $\pi_1(D-K,u)$ is an
ordered free base of $\pi_1(D-K,u)$  of  the form
$(\ell(T_1),\dots,\ell(T_n))$ where $\la T_1,\dots,T_n\ra$ is a bush in $D-K.$
\enddefinition

\medskip

\centerline{
\epsfysize=1.2in
\epsfbox{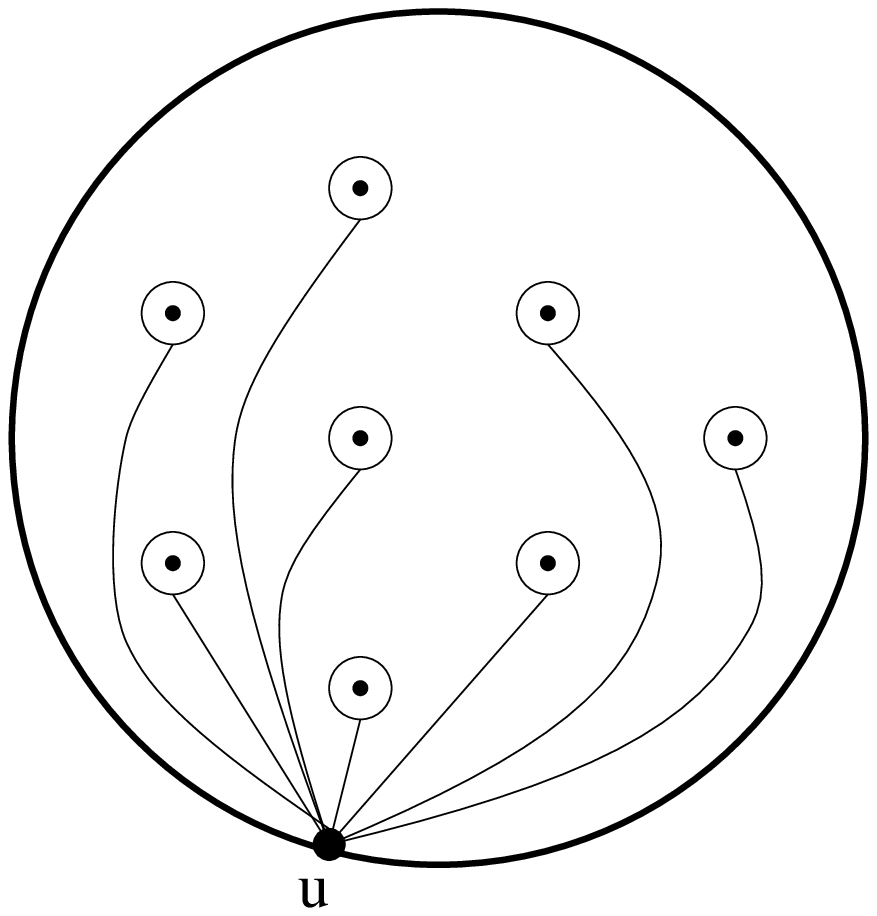}}  

\medskip

We are now able to formulate the Van Kampen theorem.
\proclaim{Van Kampen theorem} \rom{(Working format) \ \cite{VK}}
Let $\ov S\subseteq \CPt$ be a projective curve, which is transversal to the
line at infinity.
Let $S=\ov S\cap\BC^2.$
Let $\vp$ be the  braid monodromy w.r.t. $S,\pi,u.$
\quad $\vp:\pi_1(\BC-N,u)\to B_m[\BC_u^1,\BC_u^1\cap S].$
Let $\{\dl_i\}$ be a good geometric base of $\pi_1(\BC-N,u).$
Let $\{\G_j\}$ be a good geometric base of $\pi_1(\BC_u^1-S,*).$
Then:
\roster\item"(a)" $\pi_1(\BC^2-S,*)$ is generated by images of $\{\G_j\}$ in
$\pi_1(\BC^2-S,*)$ with the following relations: $\vp(\dl_i)\G_j=\G_j$\quad
$\fa i\fa j.$
\item"(b)" $\pi_1(\CPt-\ov S)\simeq
\pi_1(\BC^2-S)/\la\prod\G_j\ra.$\endroster\ep \bk
To be able to apply the Van Kampen method one has to know the actions of
$B_n[D,K]$ on $\pi_1(D-K)$ (in order to be able to compute $\vp(\dl_i)\G_j).$

One can learn how to compute the action of $B_n[D,K],$ just by considering  a simple
situation as follows.
\proclaim{Claim}
\flushpar  Assume $K=\{a,b\},$ \ $\si$ a simple path from $a$ to $b.$

\flushpar Let   $H=H(\si)\in B_n[D,K]$ be the  half-twist w.r.t. to $\si.$

\flushpar Let  $\G_a=$ a loop around $a$ counterclockwise, 
 $\G_b=$ a loop around $b$ counterclockwise. Then:
\roster\item"(a)" $(\G_a)H=\G_b.$
\item"(b)" $(\G_b)H=\G_b\G_a\G_b\1.$\endroster
\ep

\centerline{
\epsfysize=1.5cm
\epsfbox{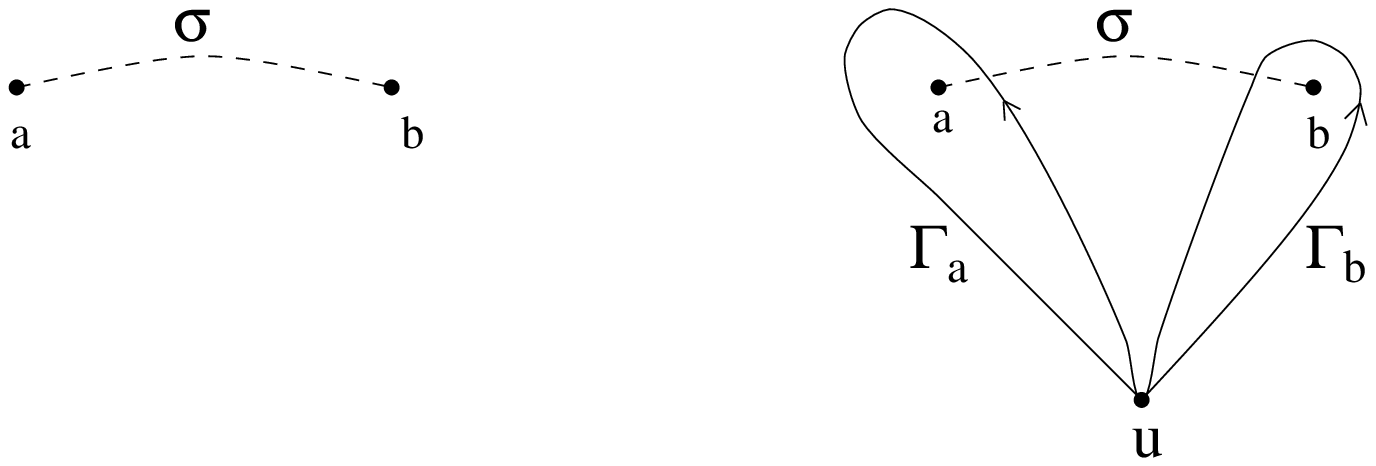}}  

\demo{Proof} 
\roster\item"(a)" Trivial.
\item"(b)" Clearly, $(\G_b)H=$

\centerline{
\epsfysize=1.5cm
\epsfbox{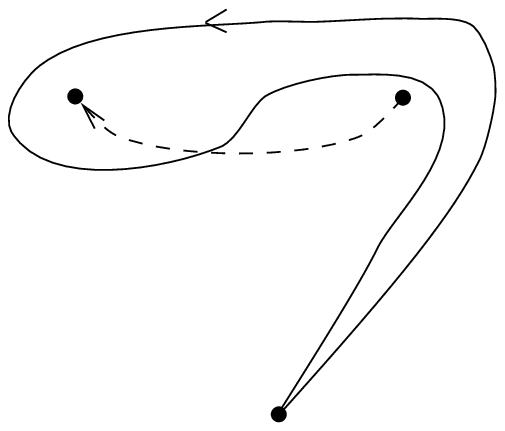}}  

Now:\newline
$\G_b\G_a=$

\centerline{
\epsfysize=1.5cm
\epsfbox{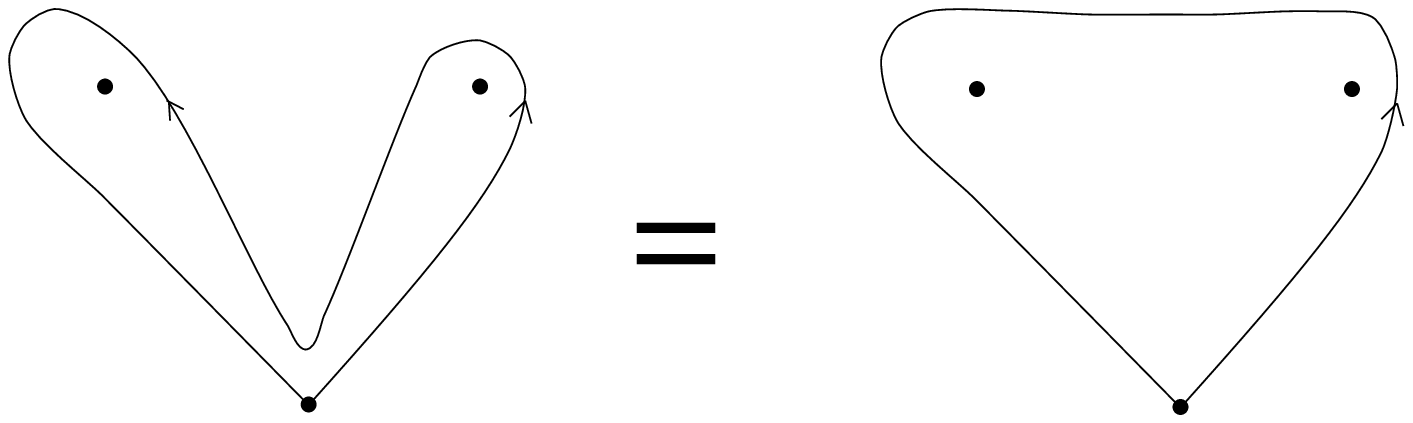}}  

\noindent
$\G_b\G_a\G_b\1=$

\centerline{
\epsfysize=1.5cm
\epsfbox{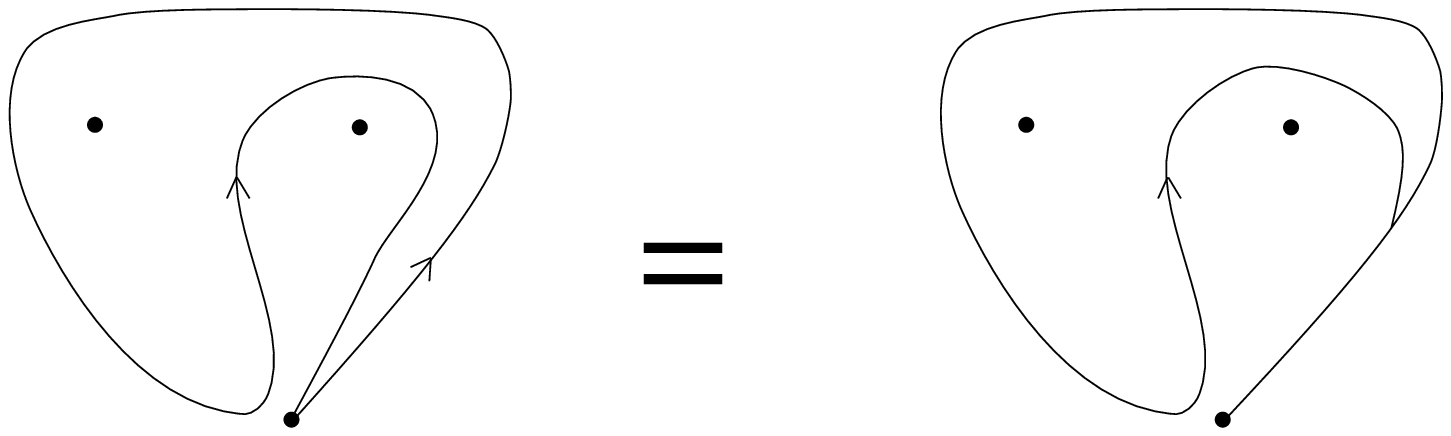}}  

$(\G_b\1$ goes {\it clockwise} around $b.$)\endroster\edm

\proclaim{Corollary} \rom{(Van Kampen for cuspidal curves)}
Let $S$ be a cuspidal curve.
The relations on $\pi_1(\BC^2-S,*))$ induced by the braid monodromy are of
the form:
$$\alignat2 &\qquad\qquad && A=B\\
&\text{or}\\
&\qquad\qquad &&AB=BA\\
&\text{or}\\
&\qquad\qquad &&ABA=BAB\endalignat$$
where $A$ and $B$ are connected to a braid $\rho(\delta)$ by the relation
$B=\rho(\delta)A$ and $A,B$ can be part of a good geometric base.
The first relation appears when $\vp(\delta)=H$, the second one when $\vp(\dl) = H^2$ and the third
one when $\vp(\dl)=H^3$\ $(H$ a half-twist). \ep

 \demo{Proof}
Let $\dl$ be an element of a $g$-base of $\pi_1(\BC^1-N).$
We want to determine the type of relation that $\vp(\dl)$ is inducing on
$\pi_1(\BC^2-S,*).$
By Van Kampen $\vp(\dl)$ induces the relations:
$$\vp(\dl)\G_j=\G_j$$
where $\{\G_j=\ell(\g_j)\}$ is a good geometric base for $\pi_1(\BC_u^1-S,*).$

Since $S$ is cuspidal, $\vp(\dl)=H^\ve$ for $\ve=1,2$ or $3.$ (See Zariski's
Theorem in the previous section.)
Thus the induced relations are
 $H^\ve(\G_j)=\G_j$\ $\forall j$.

We write $H=H(\si)$
where $\si$ connects $a$ and $b.$\edm

\demo{Case 1}    $\si$ is a straight line connecting $a$ and $b$, \
$K=\{a,b\}$ and $\g_a\si\g_b\1$ does not contain any point of $K$ in its
interior. We take $A=\G_a,$\ $B=\G_b.$

\flushpar From the previous claim we know that in $\pi_1(\BC_u^1-S,*)$:
\flushpar $H(\G_a)=\G_b$ \flushpar $H(\G_b)=\G_b\G_a\G_b\1.$ 
\flushpar$H(\G_j)=\G_j,$\quad $j\ne a,b.$
\flushpar The relation imposed on $\pi_1(\BC^2-S)$ from $\G_a=H^\ve(\G_a)$
depends on $\ve.$ \flushpar $\ve=1\Rightarrow \G_a=H(\G_a) \Rightarrow
\G_a=\G_b.$ \flushpar $\ve=2\Rightarrow
\G_a=H^2(\G_a)\Rightarrow\G_a=H(H(\G_a))=H(\G_b)=\G_b\G_a\G_b\1\Rightarrow
\G_a\G_b=\G_b\G_a$ 
\flushpar $\ve=3\Rightarrow
\G_a=H^3(\G_a)\Rightarrow
\G_a=H(H(H(\G_a))=H(\G_b\G_a\G_b\1))=\G_b\G_a\G_b\1\G_b\G_b\G_a\1\G_b\1=
\G_b\G_a\G_b\G_a\1\G_b\1\Rightarrow \G_a\G_b\G_a=\G_b\G_a\G_b$

It is easy to see  that writing $H^\ve(\G_b)=\G_b$ in $\pi_1(\BC^2-S,*)$ will
impose the same relation between $\G_a$ as $\G_b$ as the relation imposed from
$H^\ve(\G_a)=\G_a.$ The relation $H^\ve(\G_j)=\G_j$ for $j\ne a,b$ is the
trivial relation since $H(\G_j)$ equals $\G_j$ already  in
$\pi_1(\BC_u^1-S).$  Thus the realtions are of the type quoted in the lemma.
\edm
\demo{Case 2} $\si$ is any path connecting $a$ and $b$ s.t.
$K\cap\si=\{a,b\}.$

Choose a point $x$ on $\si$.
It divides $\si$ into 2 paths $\si_1$ and $\si_2.$
We choose a connection of $x$ to $*$ in $D-K.$
We call this connection $\mu(\si).$
Clearly,
$\mu(\si)\si_1\1\si\si_2\1\mu(\si)\1=(\mu(\si)\si_1)\si(\mu(\si)\si_2)\1$ has
no point  of $K$ in its interior and locally we are in the situation
of case 1.

\medskip

\centerline{
\epsfysize=1.15in
\epsfbox{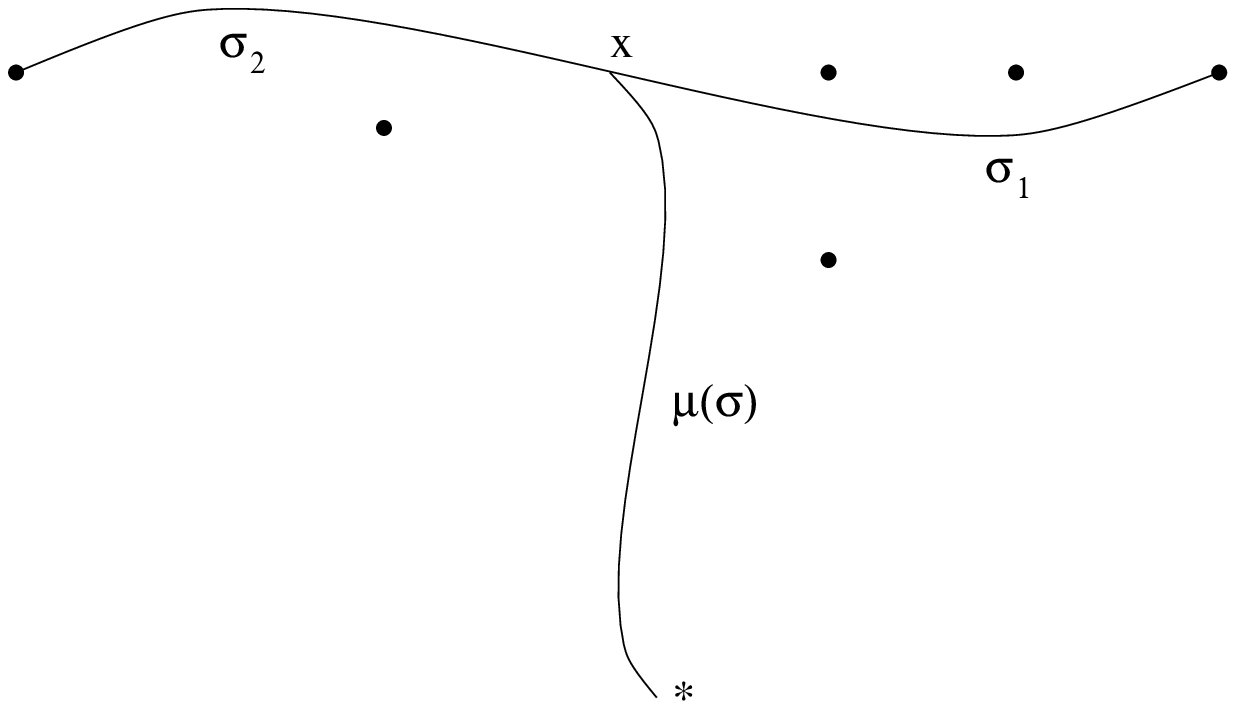}}  

\medskip

Let
$$A =\ell(\mu(\si)\si_1\1)\qquad\qquad B=\ell(\mu(\si)\si_2)$$

Since we are locally in the situation of case 1 we have $H(A)=B$ and
$H(B)=BAB\1.$  Moreover, as in  case 1
$$\align
&\ve=1\Rightarrow H(A)=A \Rightarrow   A=B\\
&\ve=2\Rightarrow H^2(A)=A\Rightarrow AB=BA\\
&\ve=3\Rightarrow H^3(A)=A\Rightarrow
ABA=BAB\endalign$$
\hfill $\square$
\edm 

\bk

\subheading{Example of simple computation of
$G=\pi(\BC^2-S,*)$ using the Van Kampen method}

\flushpar 
Let $S:y^2=x^3.$
 Clearly, $N=\{0\}$
and we take $u=1.$
$\BC_1^1\cap
S=\{-1,1\}$ and thus $\pi_1(\BC_1^1-S,u)$ is generated by $\G_1$ and
$\G_{-1}.$   $\pi_1(\BC^1-N)\simeq\la\dl\ra$ where $\dl(t)=e^{2\pi it}.$
Thus the group $\pi_1(\BC^1-N)$ induces only one relation in
$\pi_1(\BC^2-S,*).$
  In Section V we computed the braid monodromy of $\delta$ and got
 $\vp(\dl)=H^3$ where $ H=H[-1,1].$
 To
compute the relation induced on $\pi_1(\BC^1-S)$ from $\vp(\dl)$ we notice
that we are in a simple case where:
$A=\G_{-1}$\quad $B=\G_1.$
 Since $\ve=3$
the relation is $ABA=BAB.$
Thus, $ 
\pi_1(\BC^2-S)\simeq\la \G_{-1},\G_1
\ra\bigm/ ABA=BAB\simeq\la A,B\ra/ABA=BAB.$
By Artin's structure theorem we get $ \pi_1(\BC^2-S,*)\simeq
B_3.$

\remark{Remark} This example is very simple in the sense that we have only one relation while
for interesting
branch curves we have many relations ($S$ has many singular points).
In all our previous works (see \cite{MoTe2}, \cite{MoTe6}, \cite{MoTe8}, \cite{MoTe9}) we could not minimize the list of relations without first adding
more relations using invariance properties.
An invariance property of the braid monodromy is a rule with which we can replace $A$ (and $B$)
in a certain relation by a loop close to it (close enough that they almost coincide in
the degeneration).
Invariance properties are proven using the degeneration process (see, e.g., \cite{MoTe2}).

\bk

 \subheading{VIII. Some facts on the structure of $\tilde B_n$ and $\tilde B_n$-groups
(steps (e) and (f) of the algorithm)}

As pointed out in Section II, it turned out that all the new examples of $G$ and $\overline G$ are $\tilde B_n$-groups,
i.e., groups which admits an action of $\tilde B_n.$
Moreover, like $\tilde B_n$ they are extensions of  a solvable group by a symmetric one.
We shall formulate the ``almost solvability'' theorem for $\tilde B_n.$

We review first the classical Artin presentation of the braid group.
\proclaim{Theorem}
Let $X=H(\sigma_1)$ and $Y=H(\si_2)$ be \rom{2}
 half-twists.
Then
$$\si_1\cap \si_2=\emptyset\Rightarrow [X,Y]=XYX\1Y\1=1$$
$$\si_1\cap \si_2=\{\text{end point}\}\Rightarrow \la
X,Y\ra=XYXY\1 X\1Y\1=1.$$
In other words,
disjoint half-twists commute; adjacent half-twists satisfy the triple
relation.\ep
Using half-twists we build a set of generators for $B_n:$
\definition{Definition: Frame of a Braid Group}

\flushpar Take $K=\{a_1,\dots,a_n\}.$
 Let $\si_i$ be a simple path from $a_i$ to $a_{i+1}$ s.t.:
\flushpar $\si_i\cap\si_{i+1}=\{a_{i+1}\}$ and $\si_i\cap\si_j=\emptyset$ for
all $|i-j|>1.$

\flushpar  Let $X_i=H(\si_i)=$ half-twist w.r.t. $\si_i.$
  $\{X_i\}_{i=1}^{n-1}$ is called a frame of the braid group $B_n.$
\enddefinition
\remark{Remark} There is a natural epimorphism $B_n\overset\psi\to\ri S_n$ where $\psi(X_i )=
$ the transposition $(i\ i+1).$
This epimorphism ``forgets'' the diffeomorphism and only remembers the permutation of
$K.$
It is well known that a frame generation $B_n$ (Artin's theorem).
 \proclaim{Artin's
Structure Theorem}  A frame $\{X_i\}_{j=1}^n$ generates the braid group $B_n$, with only
the following relations: $$\alignat2 &X_iX_j=X_jX_i\quad &&|i-j|>1\\
&X_iX_jX_i=X_jX_iX_j\quad && |i-j|=1.\endalignat$$ \ep

We   need the following definitions for presenting the structure theorem for $\tilde B_n$.

\definition{Definitions}

$P_n=\ker(B_n\twoheadrightarrow S_n)$ where $\psi_n(X_i)=(i\ i+1)$ for some frame
$\{X_i\}$ of $B_n.$

$P_{n,0}=\ker(P_n\to Ab(B_n))$

$\tilde P_n,\tilde P_{n,0}$ the images of $P_n$ and $P_{n,0}$ in $\tilde B_n.$
\enddefinition

\proclaim{Theorem} Consider $\tilde P_n$ as a $\tilde B_n$-group.

\rom{(a)} $\tilde P_{n,0}$ is generated by a $\tilde B_n$-orbit of $\tilde X^2\tilde
Y^{-2}$ where $X$ and $Y$ are consecutive half-twists

\rom{(b)} There exist:
$$\undersetbrace \Bbb Z_2 \to {1\rightarrow (\tilde
P_n'=)}   \undersetbrace \Bbb Z^{n-1} \to { \tilde
P_{n,0}'< \tilde P} \undersetbrace \Bbb Z \to {_{n,0}'<\tilde
P}  \undersetbrace S_n \to {_n <\tilde
B_n}
$$

\rom{(c)} $\tilde P_{n,0}'=
\{1,c\}\quad c\in\Cen \tilde B_n\quad (c^2=1)$\ep
\proclaim{Corollary} $\tilde B_n$ is ``almost solvable''.
Moreover, it is an extension of a solvable group by a symmetric group.
\ep       

We
shall
not prove this theorem here. 
We only mention that the first step of the proof was the following observation.
If we have a ``good'' quadrangle in $\tilde B_n,$ i.e.

 $X_i=H(x_i)$ where $x_i$ are as above, then $\tX_1^2\tX_3^2=\tX_2^2\tX_4^2$.

\medskip

\centerline{
\epsfysize=0.95in
\epsfbox{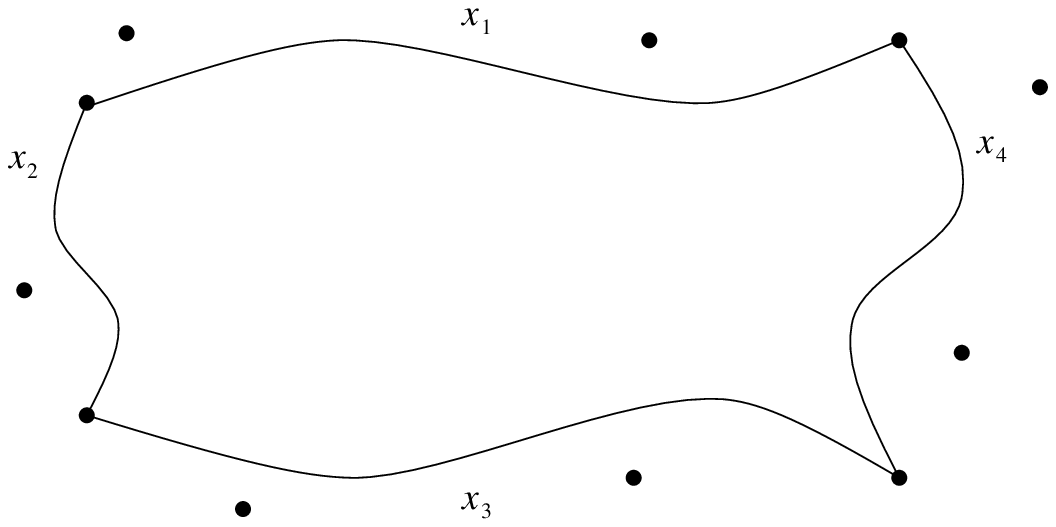}}  

\medskip

As we pointed out earlier fundamental groups turn out to be $\tilde B_n$-groups and they are also
``almost solvable''.
When studying $\tilde B_n$-groups we distinguish certain elements which we call ``prime
elements'' (i.e., $\tX^2\tY^{-2}$ in $\tilde P_n).$
Finding prime elements in a group (e.g., a fundamental  group) is the first step in 
proving that it is ``almost solvable''.

\bk

\definition {Definition:  Prime element}   
\flushpar  Let $G$ be a $\tilde B_n$-group.
We denote $b(g)$ by $g_b.$ 
An element $g\in G$ is called prime if there exists a half-twist $X\in B_n$ and $\tau\in 
\Cen(G)$
s.t. $\tau_{\tilde b}=\tau\ \forall \tilde b\in\tilde B_n,$\quad $\tau^2=1$ and
\roster\item $g_{\tilde X\1}=g\1 \tau$
\item $g_{\tX\tY\tX\1}=g\1_{\tX}g_{\tX\tY\1}\quad \forall Y$ consecutive   to $X$ 
\item $g_{\tilde Z}=g$\qquad\qquad\qquad\quad $\forall Z$ disjoint to $X$
\endroster
$X$ is called the supporting half-twist of $g.$\
\flushpar $\tau$ is called the corresponding central element.

We call these elements prime since they satisfy an existence and a uniqueness property.
We first introduce a polarization on half-twists, i.e., a direction which determines 
the beginning and end points of the path.\enddefinition

\proclaim{Existence and Uniqueness Theorem}

Let $g$ be prime supporting half-twist $X.$
Let $T$ be another half-twist. Then:
$\exists !\ h\in G$ prime and $\tilde b\in\tilde B_n$ s.t. $g_{\tilde b}=h$\quad $X_{\tilde b} 
=T$ preserving the polarization.\ep

We have proved several criteria for an element to be prime (see \cite{MoTe9} and \cite{Te1}).

\subheading{IX. The connection between fundamental groups of complements of branch curves and Galois
covers}

As pointed out in the introduction, our techniques also allow us to compute
some fundamental  groups of surfaces of general type (see
\cite{MoTe1}, \cite{MoTe2}, \cite{MoTe5}, \cite{MoRoTe}). 
These surfaces are Galois covers of generic projection to $\CPt.$ 
Recall: If
$f: X\to \CPt$ is a generic projeciton of $\deg n$ then $\tilde X$, the Galois
cover, is defined as follows:: $$\tilde
X=\overline{(X\underset    
\underbrace{\CPt\quad\CPt}_{n\
\text{times}}\to{\times\dots\times}
X)-\Dl}$$
We can compute fundamental groups of Galois covers since it can be proven
that they are quotients of a subgroups of the fundamental group of the
complement of the branch curve (see \cite{MoTe5}).  So the first steps of
computing $\pi_1(\tilde X)$ are the same as computing $\pi_1(\CPt-S).$
 In our early
attempts to find new discrete invariants for components of moduli spaces of
surfaces in the $\tau>0$ zone, we tried to use the fundamental groups of the
surface itself, believing that it can not be trivial.

In fact, until 1984 the
``Bogomolov Watershed
Conjecture'':
$$\tau>0\Rightarrow\pi_1(X)\ne
\{1\}$$
was widely believed to be true for surfaces of general type (see \cite{FH}).
In 1984 (in the process of such computations) Moishezon-Teicher disproved the
Bogomolov conjecture by constructing
counter-examples (\cite{MoTe2}).
(The surfaces we used are Galois covers of $\CP'\times\CP'.$ They  appeared in
\cite{Mi}, in which it was pointed out that these surfaces are of positive signature.) The
proof was based on computing quotients of a fundamental group of a complement of curves.
 In 1986 Chen  produced (see \cite{Ch}) 
new examples, all of
which were non-spin.
Till lately the only
known examples of spin
surfaces with
$\pi_1(X)=1,$\quad
$\tau=0$ or $\tau>0$  were the 1984
examples.
Using the Hirzebruch
surfaces we succeeded to produce infinitely many  new
examples of simply connected surfaces with $\tau>0$ and three new examples with
$\tau=0$ \cite{MoRoTe}. (Lately,  new
examples with $\tau>0$ were also produced by
Xiao, Persson and Peters
\cite{PPX}).

\bk

\newpage

\subheading{X. Galois covers of Hirzebruch surfaces: new examples}

\flushpar $F_k=$ Hirzebruch surface
of order $k$ can be schematically described as

\bigskip

\centerline{
\epsfysize=3cm
\epsfbox{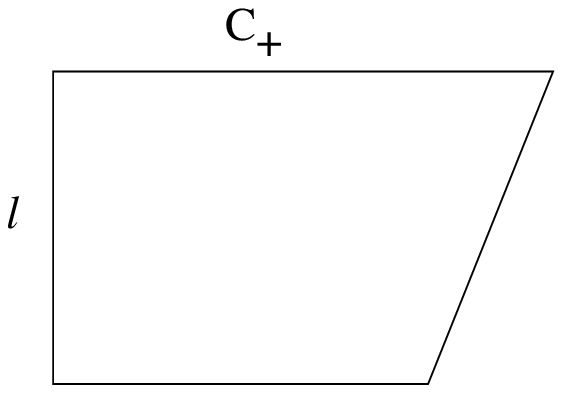}}  

\bigskip

 \flushpar where    $C_+\ell=1\quad
\ell^2=0\quad C_+^2=k.$

\flushpar 
Let $f_{ab}=f_{|a\ell+bC_{+}|}$ be the embedding of $F_k$ in $\CP^N$ w.r.t. the
full linear system
$|a\ell+bC_+|,$\quad $F_k\overset f_{ab}\to\hookrightarrow \CP^N.$
We denote $F_{k(a,b)}=f_{ab}(F_k).$
   Let $f$ be a generic projection of $F_{k(a,b)}$ to $\CPt.$
  Let $\tilde F_{k(a,b)}$ be the Galois cover of $F_{k(a,b)}$ w.r.t. $f.$\quad
 $\tilde F_{k(a,b)}$
are the new examples.  
\proclaim{Theorem}
\rom{(Moishezon, Robb, Teicher) in \cite{MoRoTe}} 
\roster\item"(a)"
For each $k$ there are infinitely many $F_{k(a,b)}$ s.t. $\tilde F_{k(a,b)}$
are simply connected surfaces of general type which are spin manifolds with
$\tau>0.$
\item"(b)"
There are $5$ surfaces which are spin with $\tau=0.$
Four of them are simply connected and the other one has  fundamental group
$\simeq \BZ_5^{48}.$ \endroster
\ep
Exact lists of $k,a,b$ and proofs can be found  in \cite{MoRoTe}.
The hardest part is computing $\pi_1(\tilde F_{k(a,b)}).$
Since 
$\pi_1(\tX)=\ker\left(\frac{\pi_1(\Bbb C^2-S)}{\la\G_j^2\ra}\to S_n\right),$ 
the first steps of computing $\pi_1(\tilde F_{k(a,b)})$ coincide with the
first steps of computing $\pi(\BC^2-S_{k(a,b)})$
where $S_{k(a,b)}$ is the branch curve of $F_{k(a,b)}\to \CPt.$
In particular, the first step is the degeneration of
the surfaces into union of planes (Step (a) of the algorithm).
We shall only present here 2 examples of degeneration.

 In fact, we present a schematic description of the degenerated object where a plane
is presented by a triangle and an intersection line between planes by an edge of a triangle. 
One can see that no 3 planes meet in a line.
The branch curve of the degenerated object is represented by the union of the edges of the
 triangles, and 
the singular points   are the intersection points of lines.
There are 2 types of singular points, depending on the number of planes/lines that come together.
In \cite{MoTe6} we describe how to determine the type of singular points of the original 
branch curve that arise 
from a singular point  of the degenerated object.
We also presented there the associated braid monodromies.
The degeneration of $V_3$ to the union of 9 planes is described in \cite{MoTe7}, the 
degeneration of $F_k(a,b)$ to the union of $2ab+kb^2$ planes is described in \cite{MoRoTe}.

\medskip

\centerline{
\epsfysize=3.8in
\epsfbox{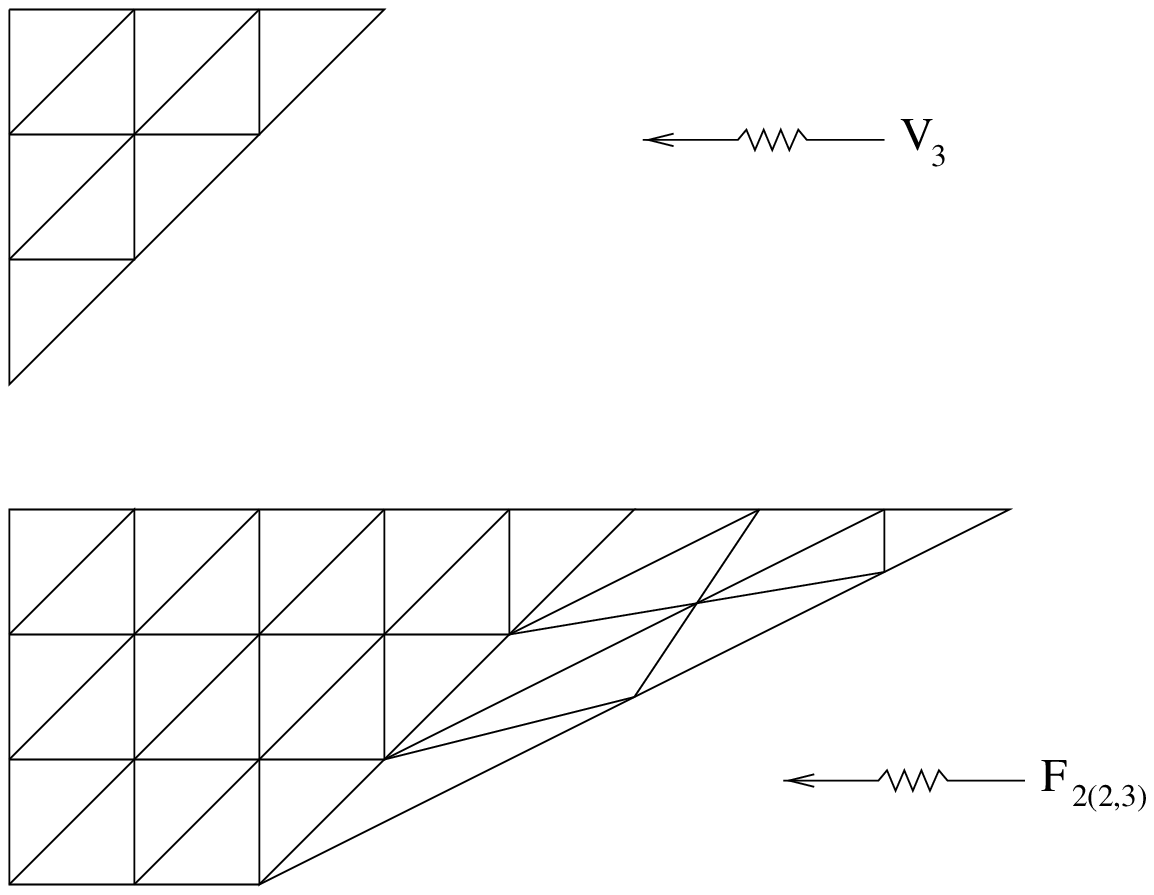}}  

\medskip

\remark{Remarks}
\roster\item $\tilde F_{0(a,b)}$ are the examples from 1984
(\cite{MoTe1},\cite{MoTe2},\cite{MoTe5}).
 \item Kotschik used $\tilde F_{0(a,b)}$
to build examples of orientation-reversing homeomorphic surfaces which are not
diffeomorphic (\cite{K}).
 \item$\tilde F_{1(a,b)}=V_b$ is the Veronese of order
$b.$ \item There is another procedure in progress to produce such examples
(\cite{MoTe11}.
\item Together with Robb and Freitag we proved lately that all
other $\tilde F_{k(a,b)}$ have finite fundamental groups which are products of
cyclic groups \cite{FRoTe}.

\newpage

\item We used \cite{MoRoTe} to produce spin surfaces of positive
signature with the same Chern classes and different fundamental groups
(see \cite{RoTe}).\endroster\endremark

 \bigskip
\tenpoint
\flushpar Department of Mathematics and Computer Science
\flushpar Bar-Ilan University
\flushpar 52900 Ramat-Gan, Israel                                                                                                                                                                                        
 \end
\definition{Definition: $\Dl_n^2$}

\flushpar Let $\{x_i\}_{i=1}^{n-1}$ be a frame of $B_n.$
\flushpar $\Dl_n^2=(x_1\cdot\dots\cdot x_{n-1})^n.$
\enddefinition

\demo{Facts}

\flushpar $\Dl_n^2\in \Cen(B_n).$
\flushpar For $n>2,$\ \ $\Dl_n^2$ generates the center.
\flushpar  For $n=2,$\ \ $B_n\simeq \BZ$ and we have: $B_n=\la x\ra$ and
$\Dl^2=x^2.$
\edm